\documentclass[12pt,oneside,reqno]{amsart}
\usepackage{graphicx}
\usepackage{mathrsfs}
\usepackage{stmaryrd}
\usepackage{amsfonts}
\usepackage{enumerate,amsmath,amssymb,amsthm}
\usepackage{booktabs} 
\usepackage{diagbox} 

\newcommand{\dif}{\mathrm{d}}

\newcommand{\be}{\begin{eqnarray}}
\newcommand{\ee}{\end{eqnarray}}
\newcommand{\ce}{\begin{eqnarray*}}
\newcommand{\de}{\end{eqnarray*}}
\newtheorem{theorem}{Theorem}[section]
\newtheorem{lemma}[theorem]{Lemma}
\newtheorem{remark}[theorem]{Remark}
\newtheorem{definition}[theorem]{Definition}
\newtheorem{proposition}[theorem]{Proposition}
\newtheorem{Examples}[theorem]{Examples}
\newtheorem{corollary}[theorem]{Corollary}

\def\t{\theta}

\def\b{\beta}
\def\d{\delta}

\def\s{\sigma}

\def\l{\lambda}

\def\[{{\Big[}}
\def\]{{\Big]}}
\def\<{{\langle}}
\def\>{{\rangle}}
\def\({{\Big(}}
\def\){{\Big)}}

\def\no{\nonumber}
\def\bt{\begin{theorem}}
\def\et{\end{theorem}}
\def\bl{\begin{lemma}}
\def\el{\end{lemma}}
\def\br{\begin{remark}}
\def\er{\end{remark}}
\def\bx{\begin{Examples}}
\def\ex{\end{Examples}}
\def\bd{\begin{definition}}
\def\ed{\end{definition}}
\def\bp{\begin{proposition}}
\def\ep{\end{proposition}}
\def\bc{\begin{corollary}}
\def\ec{\end{corollary}}

\def\cB{{\mathcal B}}

\def\cP{{\mathcal P}}

\def\mE{{\mathbb E}}

\def\mN{{\mathbb N}}

\def\mP{{\mathbb P}}

\def\mR{{\mathbb R}}

\def\mW{{\mathbb W}}

\def\sB{{\mathscr B}}
\def\sC{{\mathscr C}}

\def\sF{{\mathscr F}}

\def\geq{\geqslant}
\def\leq{\leqslant}

\begin{document}

\allowdisplaybreaks
\title{Parameter estimation of path-dependent McKean-Vlasov stochastic differential equations}

\author{Meiqi Liu$^1$ and Huijie Qiao$^{1,2}$}



\thanks{This work was supported by NSF of China (No. 11001051, 11371352) and China Scholarship Council under Grant No. 201906095034.}

\thanks{Corresponding author: Huijie Qiao, hjqiaogean@seu.edu.cn}


\date{}

\dedicatory{1. Department of Mathematics,
Southeast University,\\
Nanjing, Jiangsu 211189, P.R.China\\
2. Department of Mathematics, University of Illinois at
Urbana-Champaign\\
Urbana, IL 61801, USA}

\begin{abstract}
The work concerns a class of path-dependent McKean-Vlasov stochastic differential equations with unknown parameters. First, we prove the existence and uniqueness of these equations under non-Lipschitz conditions. Second, we construct maximum likelihood estimators of these parameters and then discuss their strong consistency. Third, a numerical simulation method for the class of path-dependent McKean-Vlasov stochastic differential equations is offered. Moreover, we estimate the errors between solutions of these equations and that of their numerical equations. Finally, we give an example to explain our result.
\end{abstract}

\maketitle \rm

\section{Introduction}

McKean-Vlasov stochastic differential equations (MVSDEs in short) are a kind of special stochastic differential equations whose coefficients depend on probability distributions of their solutions. They were first initiated by Henry P. McKean \cite{mckean} in 1966, and then were gradually studied by a lot of researchers. At present, there have been many results about MVSDEs, such as the well-posedness of the solutions in \cite{d, dq1}, the stability of strong solutions in \cite{dq2}, the well-posedness of the mild solutions and their Euler-Maruyama approximation in infinite dimension Hilbert spaces in \cite{q}, and the particle approximations method in \cite{boss}.  

As the research of MVSDEs develops, the fields of their application are becoming larger and larger. This leads to some new problems. Estimation of unknown parameters in MVSDEs is one of these problems. Now, there are many results about parameter estimation of stochastic differential equations. Let us mention some works. Liptser and Shiryayev \cite{lipster} considered the maximum likelihood estimation of It$\hat{\rm o}$ diffusions under continuous observations, while Yoshida \cite{y} estimated these diffusion processes with the maximum likelihood estimation based on discrete diffusions. In \cite{bjpn}, Bishwal obtained the exponential bound of the large deviation rate for the maximum likelihood estimator of the drift coefficients. 
Other methods of parameter estimation like martingale function estimators, nonparametric methods can be found in \cite{s,ay}.

However, because of the distributions in the drift coefficients and diffusion coefficients, the previous methods and results may not well be applied to MVSDEs. In \cite{r}, Ren and Wu proposed the least squares estimators for a class of path-dependent MVSDEs. Wen et al. \cite{w} discussed the maximum likelihood estimators on MVSDEs with the following form assuming that $\vartheta\in\mR$ is known and $\s = 1$,
 \ce
X_{t}=X_{0}+\int_{0}^{t} \int_{\mR}b\left(\t, X_{s}, y\right)\mu_s(\dif y) \dif t+\int_{0}^{t} \int_{\mR} \s\left(\vartheta, X_{s}, y\right) \mu_s(\dif y) \dif W_{s}, \quad X_{0}=x_0\in\mR,
\de
where $\t$ is a unknown parameter and $\mu_t$ is the probability distribution of $X_t$.

In this paper, we focus on the following MVSDE in a more general form 
\begin{equation}
\dif X_t = b(\t, X_{t\land\cdot}, \mu_t) \dif t + \s( X_{t\land\cdot}, \mu_t)\dif W_t, \quad X_0 =\xi,
\label{11} 
\end{equation}
where $\xi$ is a random vector. We not only construct a maximum likelihood estimator $\t_T$ for $\t$ but also prove the consistency of $\t_T$. And then we discretize Eq.(\ref{11}) and also obtain the numerical simulation of $\t_T$.

The rest of the paper is organized as follows. In Section \ref{exiuni}, we prove the existence and uniqueness of strong solutions for Eq.(\ref{11}) under non-Lipschitz conditions. The maximum likelihood estimators are constructed in Section \ref{mle}. In Section \ref{numsim}, a numerical equation of Eq.(\ref{11}) is given by interacting particles and the Euler-Maruyama method, and then the error between the MVSDE and its approximation is calculated, followed by giving a maximum likelihood estimator of the numerical equation. Finally, in Section \ref{exam}, we apply the method to a specific equation as an example, and explain our results.

The following convention will be used throughout the paper: $C$ with or without indices will denote different positive constants whose values may change from one place to another.

\section{The existence and uniqueness of path-dependent MVSDEs}\label{exiuni}

In the section, we prove the existence and uniqueness of the solutions for Eq.(\ref{11}).

Fix $T>0$. Let $C^d_T$ be the collection of all the continuous functions from $[0,T]$ to $\mR^d$. And then we equip it with the compact uniform convergence topology. Let $\cB_T^d$ be the $\sigma$-field generated by the topology. For $w\in C^d_T$, set
$$
\|w\|_{T}:=\sup_{0 \leq t \leq T}|w(t)|.
$$

Let $\sB(\mR^d)$ be the Borel $\sigma$-field on $\mR^d$. Let $\cP_2(\mR^d)$ denote the space of probability measures on $\sB(\mR^d)$ with finite second moments. That is, if $\mu\in\cP_2(\mR^d)$, then 
$$
\|\mu\|_{\lambda^2}^2:=\int_{\mR^d}(1+|x|)^2\mu(\dif x)<\infty.
$$  
And the distance of $\mu, \nu\in\cP_2(\mR^d)$ is defined as 
$$
\mW^{2}_{2}(\mu,\nu):=\inf_{\pi\in\sC(\mu_1, \mu_2)} \int_{\mR^d \times \mR^d} |x-y|^2 \pi(\dif x, \dif y),
$$
where $\sC(\mu_1, \mu_2)$ denotes the set of  all the probability measures whose marginal distributions are $\mu_1, \mu_2$, respectively. Thus, $(\cP_2(\mR^d),\mW_{2})$ is a Polish space.

Let $(\Omega,\sF,\{\sF_t\}_{0 \leq t \leq T},\mP)$ be a complete filtered probability space and $\{W_{t},t\geq0\}$ be a $m$-dimensional standard Brownian motion. Consider the following path-dependent MVSDE on $\mR^d$:
\be\left\{\begin{array}{ll}
X_t=\xi+\int_0^tb(\t, X_{s\land\cdot},\mu_s)\dif s+\int_0^t\sigma(X_{s\land\cdot},\mu_s)\dif W_s,\\
\mu_s= $the probability  distribution of\quad$X_s,
\end{array}
\label{eq1}
\right.
\ee
where $\xi$ is a $\sF_0$-measurable random vector, $\t\in\Theta\subset\mR^k$ is a unknown parameter, $b: \Theta \times C^d_T \times \cP_2(\mR^d) \mapsto \mR^d$, $\sigma: C^d_T \times \cP_2(\mR^d)\mapsto\mR^{d\times m}$ are Borel measurable. We assume:
\begin{enumerate}[($\bf{H}_1$)]
\item There exists a nonnegative constant $K_1$ such that for any $w, v\in C_T^d$, $\mu, \nu\in\cP_2(\mR^d)$

(i) 
$$
|b(\t, w,\mu) - b(\t, v,\nu) |^2+\|\sigma(w,\mu) - \s(v,\nu) \|^2 \leq K_1\left(\kappa_1(\|w-v\|_T^2) + \kappa_2\(\mW_2^2(\mu, \nu)\) \right),
$$
where $\|\cdot\|$ denotes the Hilbert-Schmidt norm of a matrix, and $\kappa_i(x), i=1, 2$ are two positive, strictly increasing, continuous concave function and satisfy $\kappa_i(0)=0$,  $\int_{0^+}\frac{1}{\kappa_1(x)+\kappa_2(x)}dx=\infty$;

(ii) 
$$
\vert b(\t, w,\mu)\vert^2+\| \s(w,\mu) \|^2 \leq K_1\left(1 + \| w\|_T^2 + \|\mu \|_{\l^2}^2 \right).
$$
\end{enumerate}

\bt\label{wellpose}
Suppose that ($\bf{H}_1$) holds and $\mE|\xi|^2<\infty$. Then Eq.(\ref{eq1}) has a unique strong solution $X$ and 
$$
\mE\left(\sup_{0 \leq t \leq T}|X(t)|^2\right)<\infty.
$$
\et
\begin{proof}
First of all, set
\be\left\{\begin{array}{ll}
X_t^{(0)}= \xi,\quad t\in[0,T],\\
X_t^{(n+1)}=\xi+\int_0^tb(\t, X_{s\land\cdot}^{(n)},\mu_s^{(n)})\dif s+\int_0^t\sigma(X_{s\land\cdot}^{(n)},\mu_s^{(n)})\dif W_s, \quad n\in\mN\cup\{0\},
\end{array}
\label{iteration}
\right.
\ee
where $\mu_s^{(n)}$ is the probability  distribution of $X_s^{(n)}$. We make use of Eq.(\ref{iteration}) to prove the well-posedness of Eq.(\ref{eq1}). 

{\bf Step 1.} We prove that the definition of Eq.(\ref{iteration}) is reasonable.

For $n=0$, $\mE\left(\sup\limits_{0 \leq t \leq T}|X_t^{(0)}|^2\right)=\mE|\xi|^2<\infty$. Assume that for $n\in\mN$,
$$
\mE\left(\sup\limits_{0 \leq t \leq T}|X_t^{(n)}|^2\right)<\infty.
$$
And then by the H$\ddot{\rm o}$lder inequality, the Burkholder-Davis-Gundy inequality  and ($\bf{H}_1$), we get that
\be
&&\mE\left(\sup_{0 \leq t \leq T}\vert X_t^{(n+1)} \vert ^2\right)\no\\
 &\leq& 3\mE|\xi|^2 +3\mE\left(\sup_{0 \leq t \leq T}\left|\int_0^tb(\t, X_{s\land\cdot}^{(n)},\mu_s^{(n)}) \dif s \right|^2\right) + 3\mE \left(\sup_{0 \leq t \leq T}\left|\int_0^t\s(X_{s\land\cdot}^{(n)},\mu_s^{(n)})\dif W_s\right|^2 \right)\no\\
&\leq &3\mE|\xi|^2 +3T \mE\int_0^T \left|b(\t,X_{s\land\cdot}^{(n)},\mu_s^{(n)})\right|^2 \dif s  + 3C\mE\int_0^T\|\s( X_{s\land\cdot}^{(n)},\mu_s^{(n)})\|^2\dif s  \no\\
&\leq &3\mE|\xi|^2 +3(T+C)K_1\mE\int_0^T\left(1 + \| X_{s\land\cdot}^{(n)}\|_T^2 + \| \mu_s^{(n)}\|_{\l^2}^2\right)\dif s\no\\
&\leq &3\mE|\xi|^2 +9(T+C)K_1T\left(1 +\mE\left(\sup\limits_{0 \leq t \leq T}|X_t^{(n)}|^2\right)\right), 
\label{linest}
\ee
where the last inequality is based on the fact that $\|\mu_s^{(n)}\|^2_{\lambda^2}\leq\mE(1+|X_s^{(n)}|)^2\leq 2\mE(1+|X_s^{(n)}|^2)$. From induction on $n$, it follows that
\ce
\mE\left(\sup\limits_{0 \leq t \leq T}|X_t^{(n)}|^2\right)<\infty, \quad n\in\mN\cup\{0\}.
\de

{\bf Step 2.} We prove the existence of the solutions to Eq.(\ref{eq1}). 

By the same deduction to that of (\ref{linest}), it holds that for $m,n\in\mN$
\be
 &&\mE\left(\sup_{0 \leq t \leq T}\vert X_t^{(n+1)} - X_t^{(m+1)} \vert ^2 \right)\no\\
 &\leq& 2T\mE \int_0^T \vert b(\t, X_{s\land\cdot}^{(n)}, \mu_s^{(n)}) -b(\t, X_{s\land\cdot}^{(m)}, \mu_s^{(m)})\vert^2 \dif s\no\\
 &&+ 2C\mE\int_0^T \|\s(X_{s\land\cdot}^{(n)}, \mu_s^{(n)}) -\s( X_{s\land\cdot}^{(m)}, \mu_s^{(m)} )\|^2 \dif s \no\\
&\leq& 2(T+C)K_1\mE \int_0^T\bigg(\kappa_1(\| X_{s\land\cdot}^{(n)} - X_{s\land\cdot}^{(m)} \|_T^2)\no\\
&&\quad +\kappa_2(\mW_2^2(\mu_s^{(n)},\mu_s^{(m)}))\bigg)\dif s\no\\
&\leq& 2(T+C)K_1\int_0^T\bigg[\kappa_1\left(\mE\left(\sup\limits_{0 \leq u\leq s}|X_u^{(n)}-X_u^{(m)}|^2\right)\right)\no\\
&&+\kappa_2\left(\mE\left(\sup\limits_{0 \leq u\leq s}|X_u^{(n)}-X_u^{(m)}|^2\right)\right)\bigg]\dif s,
\label{distesti}
\ee 
where the last step is based on the Jensen inequality and the fact that 
$$
\mW_2^2(\mu_s^{(n)},\mu_s^{(m)}) \leq \mE|X_s^{(n)}-X_s^{(m)}|^2 \leq \mE\left(\sup\limits_{0 \leq u \leq s}|X_u^{(n)}-X_u^{(m)}|^2\right). 
$$
Set
$$
g(t):=\lim\limits_{n,m\rightarrow\infty}\mE\left(\sup\limits_{0 \leq u\leq t}|X_u^{(n)}-X_u^{(m)}|^2\right),
$$
and then (\ref{distesti}) admits us to have that
$$
g(T)\leq 2(T+C)K_1\int_0^T\(\kappa_1(g(s))+\kappa_2(g(s))\)\dif s.
$$
Thus, by \cite[Lemma 3.6]{dq1}, one can get $g(T)=0$. That is, $\{X^{(n)}\}$ is a Cauchy sequence in the space $L^2(\Omega,\sF,\mP, C_T^d)$. From this, we know that there exists a $X\in L^2(\Omega,\sF,\mP, C_T^d)$ such that 
\be
\lim\limits_{n\to\infty}\mE\left(\sup_{0 \leq t \leq T}\vert X_t^{(n)} - X_t\vert ^2\right)=0.
\label{limit01}
\ee
Note that 
$$
\sup_{0 \leq t \leq T}\mW_2^2(\mu_t^{(n)},\mu_t) \leq \sup_{0 \leq t \leq T}\mE|X_t^{(n)}-X_t|^2 \leq \mE\left(\sup\limits_{0 \leq t \leq T}|X_t^{(n)}-X_t|^2\right). 
$$
So, we conclude that
\be
\lim\limits_{n\to\infty}\sup_{0 \leq t \leq T}\mW_2^2(\mu_t^{(n)},\mu_t)=0.
\label{limit02}
\ee
And then (\ref{limit01})-(\ref{limit02}) imply that for $\forall t \in [0,T]$,
\ce
&&\int_0^tb(\t, X_{s\land\cdot}^{(n)}, \mu_s^{(n)}) \dif s  \to \int_0^t b(\t, X_{s\land\cdot}, \mu_s) \dif s ,\quad a.s.,\\
&&\int_0^t \s(X_{s\land\cdot}^{(n)}, \mu_s^{(n)})\dif W_s  \to \int_0^t\s( X_{s\land\cdot}, \mu_s) \dif W_s ~\mbox{in}~ L^2(\Omega, \sF_t,\mP).
\de
Therefore, taking the limit on two hand sides of Eq.(\ref{iteration}) as $n \to \infty$, we have that 
$$
X_t=\xi+\int_0^tb(\t, X_{s\land\cdot},\mu_s)\dif s+\int_0^t\sigma(X_{s\land\cdot},\mu_s)\dif W_s,
$$
that is, $X$ is a solution of Eq.(\ref{eq1}). 

{\bf Step 3.} We prove the uniqueness of the solutions to Eq.(\ref{eq1}). 

Suppose that $X$ and $\hat{X}$ are two solutions to Eq.(\ref{eq1}). And then by the similar calculation to that of (\ref{distesti}), it holds that
\ce
 \mE\left(\sup_{0 \leq t \leq T}\vert X_t - \hat X_t \vert ^2\right) 
&\leq& 2T\mE \int_0^T\vert b(\t, X_{s\land\cdot}, \mu_s) -b(\t, \hat{X}_{s\land\cdot}, \hat \mu_s)\vert^2 \dif s   \\
&&+ \ 2C\mE\int_0^T \|\s(X_{s\land\cdot}, \mu_s) -\s( \hat{X}_{s\land\cdot}, \hat \mu_s）)\|^2 \dif s \\
&\leq& 2(T+C)K_1\mE\int_0^T\(\kappa_1(\| X_{s\land\cdot} - \hat{X}_{s\land\cdot}\|_T^2) + \kappa_2(\mW_2^2(\mu_s,\hat \mu_s))\) \dif s\\
&\leq&2(T+C)K_1\int_0^T \bigg(\kappa_1\left(\mE\left(\sup_{0 \leq u \leq s}\vert X_u - \hat X_u \vert ^2\right)\right)\\
&&+\kappa_2\left(\mE\left(\sup_{0 \leq u \leq s}\vert X_u - \hat X_u \vert ^2\right)\right)\bigg)\dif s,
\de 
which together with \cite[Lemma 3.6]{dq1} yields that
$$
\mE\left(\sup_{0 \leq t \leq T}\vert X_t - \hat X_t \vert ^2\right)=0.
$$
That is, $X_t=\hat X_t$ for all $t\in[0,T]$ and almost all $\omega$. The proof is complete.
\end{proof}

\section{The maximum likelihood estimation of path-dependent MVSDEs}\label{mle}

In the section, we assume ($\bf{H}_1$) and $d=m=k=1$. And then Eq.(\ref{eq1}) has a unique solution $X^{\t}$. We construct a maximum likelihood estimator of $\t$ and prove its properties. $C_T:=C_T^1$.

Assume:
\begin{enumerate}[($\bf{H}_2$)] 
\item For any $w\in C_T, \mu\in\cP_2(\mR)$, $\s(w,\mu)\neq 0$ and  
$$
\left| \frac{b(\t,w,\mu)}{\s(w,\mu)}\right|\leq K_2,
$$
where $K_2\geq 0$ is a constant.
\end{enumerate}

Let $\t_0$ be the true value of $\t$. Let $\mP_{\t}^{T}, \mP_{\t_0}^{T}$ be the distributions of $(X^\t_t)_{t\in[0,T]}$ and $(X^{\t_0}_t)_{t\in[0,T]}$, respectively. Thus, under ($\bf{H}_2$), it follows from \cite[Theorem 7.19, P. 294]{lipster} that $\mP_{\theta}^{T} \ll \mP_{\theta_{0}}^{T}$. Define a maximum likelihood function of $\t$ as
\ce
L_{T}(\theta) &:=&\frac{\dif {\mP_{\theta}^{T}}}{\dif \mP^{T}_{\theta_{0}}} \\
&=&\exp\left\{{\int_{0}^{T} \frac{1}{\s^2(X_{t\land\cdot}^{\theta_{0}}, \mu^{\theta_{0}}_{t}) }\(b(\t, X_{t\land\cdot}^{\theta_{0}}, \mu^{\theta_{0}}_t) - b(\t_0, X_{t\land\cdot}^{\theta_{0}}, \mu^{\theta_{0}}_t) \)\dif X_{t}^{\theta_0}}\right.\\
&&\qquad\qquad -\left.{\frac{1}{2}\int_{0}^{T} \frac{1}{\s^2(X_{t\land\cdot}^{\theta_{0}}, \mu^{\theta_{0}}_{t}) }\(b^2(\t, X_{t\land\cdot}^{\theta_{0}}, \mu^{\theta_{0}}_t) - b^2(\t_0, X_{t\land\cdot}^{\theta_{0}}, \mu^{\theta_{0}}_t) \)\dif t}\right\}\\
&=&\exp\left\{{\int_{0}^{T} \frac{1}{\s(X_{t\land\cdot}^{\theta_{0}}, \mu^{\theta_{0}}_{t}) }\(b(\t, X_{t\land\cdot}^{\theta_{0}}, \mu^{\theta_{0}}_t) - b(\t_0, X_{t\land\cdot}^{\theta_{0}}, \mu^{\theta_{0}}_t) \)\dif W_{t} }\right.\\
&&\qquad\qquad-\left.{\frac{1}{2}\int_{0}^{T} \frac{1}{\s^2(X_{t\land\cdot}^{\theta_{0}}, \mu^{\theta_{0}}_{t}) }\(b(\t, X_{t\land\cdot}^{\theta_{0}}, \mu^{\theta_{0}}_t) - b(\t_0, X_{t\land\cdot}^{\theta_{0}}, \mu^{\theta_{0}}_t) \)^2\dif t}\right\},
\de
where $\mu^\t_t, \mu^{\t_0}_t$ are the distributions of $X^\t_t, X^{\t_0}_t$, respectively. So, the maximum likelihood estimator of $\t$ is given by
\ce
\t_T := \arg\max_{\t \in \Theta}{L_{T}(\t)}.
\de

Next, we study some properties of the maximum likelihood estimator $\t_T$. To do this, we assume more:
\begin{enumerate}[($\bf{H}_3$)] 
\item For any $w\in C_T, \mu\in\cP_2(\mR)$, $b(\t, w, \mu)$ is one-to-one and continuous in $\t$.
\end{enumerate}

\bt\label{et2}(The strong consistency)\\
Under the assumptions $(\bf{H}_1)$-$(\bf{H}_3)$, it holds that
$$
\t_T\stackrel{a.s.}{\longrightarrow}\t_0,  \quad T\to \infty.
$$
\et
\begin{proof}
Set 
$$
l_{T}(\t) := \log L_T(\t) =\displaystyle\log\frac{\dif {\mP_{\t}^{T}}}{\dif \mP_{\t_{0}}^{T}}.
$$
And then it holds that for $\d > 0$,
\ce
l_{T}(\t_0+\d) - l_T(\t_0) &=&\log\frac{\dif {\mP_{\t_0+\d}^{T}}}{\dif \mP_{\t_{0}}^{T}} \\
&=&\int_{0}^{T} \frac{1}{\s(X_{t\land\cdot}^{\theta_{0}}, \mu^{\theta_{0}}_{t}) }\(b(\t_0+\d, X_{t\land\cdot}^{\theta_{0}}, \mu^{\theta_{0}}_t) - b(\t_0, X_{t\land\cdot}^{\theta_{0}}, \mu^{\theta_{0}}_t) \)\dif W_{t} \\
&&-\frac{1}{2}\int_{0}^{T} \frac{1}{\s^2(X_{t\land\cdot}^{\theta_{0}}, \mu^{\theta_{0}}_{t}) }\(b(\t_0+\d, X_{t\land\cdot}^{\theta_{0}}, \mu^{\theta_{0}}_t) - b(\t_0, X_{t\land\cdot}^{\theta_{0}}, \mu^{\theta_{0}}_t) \)^2\dif t\\
&=:&{\int_{0}^{T} \Gamma_t^{\t_0} \dif W_{t} }-\frac{1}{2}\int_{0}^{T} (\Gamma_t^{\t_0})^2\dif t, 
\de
where 
$$
\Gamma_t^{\t_0}  := \frac{1}{\s(X_{t\land\cdot}^{\theta_{0}}, \mu^{\theta_{0}}_{t}) }\(b(\t_0+\d, X_{t\land\cdot}^{\theta_{0}}, \mu^{\theta_{0}}_t) - b(\t_0, X_{t\land\cdot}^{\theta_{0}}, \mu^{\theta_{0}}_t) \). 
$$
Note that 
$$
\left[\int_{0}^{\cdot} \Gamma_t^{\t_0} \dif W_{t}\right]_T=\int_{0}^{T}\left|\Gamma_t^{\t_0}\right|^2\dif t,
$$ 
where $[\cdot]$ stands for the quadratic variation of $\cdot$. Thus, by the time change, we know that 
$$
\tilde{W}_t:=\int_{0}^{A_t} \Gamma_s^{\t_0} \dif W_{s}
$$
is a $(\sF_{A_t})_{t\geq 0}$-adapted Brownian motion, where $A_t$ is the inverse function of $\int_{0}^{t}\left|\Gamma_s^{\t_0}\right|^2\dif s$. So,
\be
\frac{l_T(\t_0+\delta)-l_T(\t_0)}{\int_{0}^{T} (\Gamma_t^{\t_0})^2\dif t}
=\frac{\int_{0}^{T} \Gamma_t^{\t_0}\dif W_t }{\int_{0}^{T} (\Gamma_t^{\t_0})^2\dif t}-\frac{1}{2}=\frac{\tilde{W}_{A_T^{-1}}}{A_T^{-1}}-\frac{1}{2} \stackrel{a.s.}{\longrightarrow}-\frac{1}{2}, \quad T\rightarrow\infty,
\label{limit1}
\ee
where the last step is based on the strong law of large numbers for Brownian motions. By the same deduction to that of (\ref{limit1}), one can get that
\be
\frac{l_T(\t_0-\delta)-l_T(\t_0)}{\int_{0}^{T} (\Gamma_t^{\t_0})^2\dif t}\stackrel{a.s.}{\longrightarrow}-\frac{1}{2}, \quad T\rightarrow\infty.
\label{limit2}
\ee
Combining (\ref{limit1}) with (\ref{limit2}), we obtain that
\be
\frac{l_T(\t_0\pm\delta)-l_T(\t_0)}{\int_{0}^{T} (\Gamma_t^{\t_0})^2\dif t}\stackrel{a.s.}{\longrightarrow}-\frac{1}{2}, \quad T\rightarrow\infty.
\label{limit}
\ee

Next, we observe (\ref{limit}). It follows from (\ref{limit}) that for $\delta$ and $\theta_0$, there exists some $t_0>0$ such that
\begin{equation}
l_{T}(\theta_{0}\pm\delta) <l_{T}(\theta_{0}), \quad T \geq t_0, \quad a.s..
\label{26} 
\end{equation}
Besides, by ($\bf{H}_3$), we know that $l_T(\theta)$ is continuous on $[\t_0 - \d, \t_0 + \d]$. So, there exists a $\t^*\in[\t_0 - \d, \t_0 + \d]$ such that $l_T(\t^*)$ is the maximum value of $l_T(\theta)$ on $[\t_0 - \d, \t_0 + \d]$. That is, $\t_T=\t^*$ for $\Theta=[\t_0 - \d, \t_0 + \d]$. Based on (\ref{26}), it holds that $\t_T \neq \t_0 \pm \d $ for $T\geq t_0$. Thus, $\theta_T \to \theta_0$ as $T\to\infty$. The proof is over.
\end{proof}

\section{The numerical simulation of path-dependent MVSDEs}\label{numsim}

In the section, we introduce the numerical simulation of Eq.(\ref{eq1}) under ($\bf{H}_1$) and estimate the error between the solution of Eq.(\ref{eq1}) and that of the numerical equation under Lipschitz conditions.

First of all, for $N\in\mN$ consider these following MVSDEs
\be\left\{\begin{array}{ll}
\dif X_{t}^{i,N}= b\left(\t,  X_{t\land\cdot}^{i,N}, \mu_t^N\right) \dif t +\sigma\left( X_{t\land\cdot}^{i,N}, \mu_t^N \right) \dif W_t^i ,\\
X_{0}^{i,N}=\xi, \quad i = 1, 2,\ldots, N,
\end{array}
\label{31} 
\right.
\ee
where $\mu_t^N:=\dfrac{1}{N}\sum\limits_{j=1}^{N}\d_{X_t^{j,N}}$, $\d_{X_t^{j,N}}$ is the Dirac measure at $X_t^{j,N}$, and $W_t^i, i = 1, 2,\ldots, N$ are $N$ mutually independent $m$-dimensional standard Brownian motions. By Theorem \ref{wellpose}, under ($\bf{H}_1$) we know that Eq.(\ref{31}) has a unique solution $X_{t}^{i,N}$. And then we construct  the following numerical simulation equation: for $M\in\mN$
\be\left\{\begin{array}{ll}
Y_{0}^{i}=\xi,\\
Y_t^{i}=Y_{t_k}^{i}+b\left(\theta,  Y_{t_k\land\cdot}^{i}, \mu_{t_k}^M\right)(t-t_k) +\sigma\left( Y_{t_k\land\cdot}^{i}, \mu_{t_k}^M\right)(W^i_t-W^i_{t_k}), \quad t\in[t_k, t_{k+1}],
\end{array}
\label{34} 
\right.
\ee
where $t_k:=k\frac{T}{M}$, $\mu_{t_k}^M :=\dfrac{1}{N}\sum\limits_{j=1}^{N}\delta_{Y_{t_k}^j} $ for $k= 0,...,M-1$. In order to estimate the error between the solution of Eq.(\ref{34}) and the solution of Eq.(\ref{eq1}), we also introduce the following MVSDE:
\be
X_t^i = \xi+\int_0^t b(\t, X_{s\land\cdot}^i, \mu^i_s) \dif s + \int_0^t\s( X_{s\land\cdot}^i, \mu^i_s)\dif W_s^i,
\label{33} 
\ee
where $\mu^i_s$ is the distribution of $X_s^i$. Note that the solution of Eq.(\ref{33}) has the same distribution to that of the solution for Eq.(\ref{eq1}). Therefore, we compute the distance between $X_t^i$ and $Y_t^{i}$ to estimate the error between $X_t$ and $Y_t^{i}$. To do this, we need stronger assumptions than ($\bf{H}_1$). Assume:
\begin{enumerate}[($\bf{H}^{\prime}_1$)]
\item There exists a nonnegative constant $K'_1$ such that for any $w, v\in C_T^d$, $\mu, \nu\in\cP_2(\mR^d)$

(i) 
$$
|b(\t, w,\mu) - b(\t, v,\nu) |^2+\|\sigma(w,\mu) - \s(v,\nu) \|^2 \leq K'_1\left( \|w-v\|^2 + \mW_2^2(\mu, \nu) \right),
$$

(ii) 
$$
\vert b(\t, w,\mu)\vert^2+\| \s(w,\mu) \|^2 \leq K'_1\left(1 + \| w\|^2 + \|\mu \|_{\l^2}^2 \right).
$$
\end{enumerate}

\bt
Suppose that ($\bf{H}'_1$) holds and $\mE|\xi|^p<\infty$ for $p>4$. Then it follows that
\be
\sup_{1 \leq i \leq N} \mE\left[ \sup_{0 \leq t \leq T}\vert X_t^i - Y_t^i \vert ^2\right] \leq C\Gamma_N+C \frac{T}{M}\left(\frac{T}{M}+C\right),
\label{errest}
\ee
where the constant $C>0$ is independent of $N, M$ and \ce
\Gamma_N:=\left\{\begin{array}{ll}
N^{-1/2}, \qquad\qquad d<4,\\
N^{-1/2}\log N, ~\quad d=4,\\
N^{-1/d}, \qquad\qquad d>4.
\end{array}
\right.
\de
\et
\begin{proof}
Note that
\be
\sup_{1 \leq i \leq N} \mE\left[ \sup_{0 \leq t \leq T}\vert X_t^i - Y_t^i \vert ^2\right] &\leq&2\sup_{1 \leq i \leq N} \mE\left[ \sup_{0 \leq t \leq T}\vert X_t^i - X_t^{i,N} \vert ^2\right] \no\\
&&+2\sup_{1 \leq i \leq N} \mE\left[ \sup_{0 \leq t \leq T}\vert X_t^{i,N} - Y_t^i \vert ^2\right]\no\\
&=:& I_1+I_2.
\label{35} 
\ee
For $I_1$, it follows from the same deduction as that of (\ref{linest}) that for $\forall i=1,\ldots, N,$
\ce
\mE\left(\sup_{0 \leq t\leq T} \vert X_t^i - X_t^{i,N} \vert ^2\right)&\leq& 2\mE\sup_{0 \leq t\cdot \leq T}t\int_0^t\left|b(\t, X_{s\land\cdot}^i, \mu^i_s) - b(\t, X_{s\land\cdot}^{i,N}, \mu_s^N)\right|^2 \dif s \\
&&+2C\mE\int_0^T \left\|\s(X_{s\land\cdot}^i, \mu^i_s) -\s( X_{s\land\cdot}^{i,N}, \mu_s^N)\right\|^2 \dif s \\
&\leq& 2(T+C)K'_1\mE \int_0^T\(\|X_{s\land\cdot}^i -X_{s\land\cdot}^{i,N}\|_T^2+\mW_2^2(\mu^i_s,\mu_s^N)\) \dif s \\
&\leq& 2(T+C)K'_1\int_0^T\mE\left(\sup_{0 \leq r \leq s}\vert X_r^i - X_r^{i,N} \vert ^2\right)\dif s \\
&&+2(T+C)K'_1\int_0^T\mE\left(\mW_2^2(\mu^i_s,\mu_s^N)\right) \dif s.
\de
Gronwall's inequality admits us to obtain that
\ce
\mE\left(\sup_{0 \leq t\leq T} \vert X_t^i - X_t^{i,N} \vert ^2\right)&\leq&2(T+C)K'_1\int_0^T  \mE \mW_2^2(\mu^i_s,\mu_s^N) \dif s \cdot \exp\{2(T+C)K'_1T\} \\
&\leq&C\Gamma_N,
\de
where the last inequality is based on \cite[Theorem 5.8, P. 362]{cf}, and furthermore
\be
I_1\leq C\Gamma_N.
\label{i1}
\ee

For $I_2$, by the similar deduction to that of (\ref{linest}), it holds that
\ce
\mE\left[ \sup_{0 \leq t \leq T}\vert X_t^{i,N} - Y_t^i \vert ^2\right]
&\leq&2T\int_0^T\mE\left|b(\t,X_{s\land\cdot}^{i,N},\mu_s^N)-b(\t,Y^i_{\eta(s)\land\cdot},\mu_{\eta(s)}^M)\right|^2\dif s\\
&&+2C\int_0^T\mE\left|\s(X_{s\land\cdot}^{i,N},\mu_s^N)-\s(Y^i_{\eta(s)\land\cdot},\mu_{\eta(s)}^M)\right|^2\dif s\\
&\leq&(2T+2C)K'_1\int_0^T\mE\left(\|X_{s\land\cdot}^{i,N}-Y^i_{\eta(s)\land\cdot}\|_T^2+\mW_2^2(\mu_s^N,\mu_{\eta(s)}^M)\right)\dif s\\
&\leq&(2T+2C)K'_1\int_0^T\mE\(2\|X_{s\land\cdot}^{i,N}-Y_{s\land\cdot}^i\|_T^2+2\|Y_{s\land\cdot}^i-Y^i_{\eta(s)\land\cdot}\|_T^2\\
&&+2\mW_2^2(\mu_s^N,\mu^M_s)+2\mW_2^2(\mu^M_s,\mu_{\eta(s)}^M)\)\dif s\\
&\leq&8(T+C)K'_1\int_0^T\mE\(\sup\limits_{0\leq r\leq s}|X_r^{i,N}-Y_r^i|^2\)\dif s\\&&+8(T+C)K'_1T\sup\limits_k\mE\left(\sup\limits_{t_k\leq r\leq t_{k+1}}|Y^i_r-Y^i_{t_k}|^2\right),
\de
where $\eta(s)=t_k, s\in[t_k, t_{k+1}]$ and the following fact is used:
$$
\mE\mW_2^2(\mu^N_s,\mu_s^M)\leq\mE\left(\dfrac{1}{N}\sum\limits_{j=1}^{N}|X_s^{j,N}-Y_s^j|^2\right)=\mE|X_s^{i,N}-Y_s^i|^2.
$$
The Gronwall inequality admits us to obtain that
\be
\mE\left[ \sup_{0 \leq t \leq T}\vert X_t^{i,N} - Y_t^i \vert ^2\right]\leq 8(T+C)K'_1T\sup\limits_k\mE\left(\sup\limits_{t_k\leq r\leq t_{k+1}}|Y^i_r-Y^i_{t_k}|^2\right)e^{8(T+C)K'_1T}.
\label{mides}
\ee

In the following, we estimate $\mE\left(\sup\limits_{t_k\leq r\leq t_{k+1}}|Y^i_r-Y^i_{t_k}|^2\right)$. By (\ref{34}), it holds that
\be
\mE\left(\sup\limits_{t_k\leq r\leq t_{k+1}}|Y^i_r-Y^i_{t_k}|^2\right)&\leq&2\mE\left(\sup\limits_{t_k\leq r\leq t_{k+1}}\left|\int_{t_k}^r b\left(\theta,  Y_{t_k\land\cdot}^{i}, \mu_{t_k}^M\right)\dif u \right|^2\right)\no\\
 &&+2\mE\left(\sup\limits_{t_k\leq r\leq t_{k+1}}\left|\int_{t_k}^r\sigma\left( Y_{t_k\land\cdot}^{i}, \mu_{t_k}^M\right)\dif W^i_u\right|^2\right)\no\\
&\leq&2\frac{T^2}{M^2}\mE|b\left(\theta,  Y_{t_k\land\cdot}^{i}, \mu_{t_k}^M\right)|^2
+2C\int_{t_k}^{t_{k+1}}\mE\|\sigma\left( Y_{t_k\land\cdot}^{i}, \mu_{t_k}^M\right)\|^2\dif u\no\\
&\leq&2\frac{T^2}{M^2}\mE|b\left(\theta,  Y_{t_k\land\cdot}^{i}, \mu_{t_k}^M\right)|^2+2C\frac{T}{M}\mE\|\sigma\left( Y_{t_k\land\cdot}^{i}, \mu_{t_k}^M\right)\|^2\no\\
&\leq&2\frac{T}{M}\left(\frac{T}{M}+C\right)\mE\left(1+\|Y_{t_k\land\cdot}^{i}\|_T^2+\|\mu_{t_k}^M\|^2_{\lambda^2}\right)\no\\
&\leq&2\frac{T}{M}\left(\frac{T}{M}+C\right)\left(1+\mE\left(\sup\limits_{0\leq r\leq t_k}|Y^i_r|^2\right)+\mE\(1+|Y^i_{t_k}|\)^2\right)\no\\
&\leq&6\frac{T}{M}\left(\frac{T}{M}+C\right)\left(1+\mE\left(\sup\limits_{0\leq r\leq t_k}|Y^i_r|^2\right)\right),
\label{ttkes}
\ee
where in the last second inequality we use the fact that
$$
\mE\|\mu_{t_k}^M\|^2_{\lambda^2}=\frac{1}{N}\sum\limits_{j=1}^{N}\mE\int_{\mR^d}(1+|x|)^2\delta_{Y_{t_k}^j} (\dif x)=\frac{1}{N}\sum\limits_{j=1}^{N}\mE(1+|Y_{t_k}^j|)^2=\mE\(1+|Y^i_{t_k}|\)^2.
$$
Besides, from the similar deduction to that of (\ref{linest}), it follows that
\ce
&&\mE\left(\sup\limits_{0\leq t\leq T}|Y^i_t|^2\right)\\
&\leq&3\mE|\xi|^2+3T\mE\int_0^T\left|b(\t,Y^i_{\eta(s)\land\cdot},\mu_{\eta(s)}^M)\right|^2ds+3\mE\int_0^T\left\|\sigma(Y^i_{\eta(s)\land\cdot},\mu_{\eta(s)}^M)\right\|^2ds\no\\
&\leq&3\mE|\xi|^2+3(T+1)\mE\int_0^TK'_1\left(1+\|Y^i_{\eta(s)\land\cdot}\|_T^2+\|\mu_{\eta(s)}^M\|^2_{\lambda^2}\right)\dif s\no\\
&\leq&3\mE|\xi|^2+3(T+1)\int_0^TK'_1\left(1+\mE\left(\sup\limits_{0\leq u\leq s}|Y^i_u|^2\right)+2\mE(1+|Y^i_{\eta(s)}|^2)\right)\dif s\no\\
&\leq&3\mE|\xi|^2+9(T+1)TK'_1+9(T+1)K'_1\int_0^T\mE\left(\sup\limits_{0\leq u\leq s}|Y^i_u|^2\right)\dif s.
\de
The Gronwall inequality admits us to obtain that 
\be
\mE\left(\sup\limits_{0\leq t\leq T}|Y^i_t|^2\right)\leq C.
\label{2semo}
\ee
Combing (\ref{mides})-(\ref{2semo}), we have that 
\ce
\mE\left[ \sup_{0 \leq t \leq T}\vert X_t^{i,N} - Y_t^i \vert ^2\right]\leq C\frac{T}{M}\left(\frac{T}{M}+C\right),
\de
and furthermore
\be
I_2\leq C\frac{T}{M}\left(\frac{T}{M}+C\right).
\label{i2es}
\ee

Finally, from (\ref{35}) (\ref{i1}) (\ref{i2es}), it follows that
\ce
\sup_{1 \leq i \leq N} \mE\left[ \sup_{0 \leq t \leq T}\vert X_t^i - Y_t^i \vert ^2\right] \leq C \Gamma_N+C \frac{T}{M}\left(\frac{T}{M}+C\right).
\de
The proof is complete.
\end{proof}

\bigskip

Next, we construct a maximum likelihood estimator of the parameter $\t$. Let $d=m=k=1$. Assume that ($\bf{H}'_1$)-($\bf{H}_2$) hold. And then define the maximum likelihood function 
\ce
L^M_{T}(\theta) 
&:=&\exp\left\{{\int_{0}^{T} \frac{1}{\s(Y_{\eta(t)\land\cdot}^i, \mu_{\eta(t)}^M) }\left(b(\t, Y_{\eta(t)\land\cdot}^i, \mu_{\eta(t)}^M) - b(\t_0, Y_{\eta(t)\land\cdot}^i, \mu_{\eta(t)}^M) \right)\dif W^i_{t} }\right.\no\\
&&\qquad -\left.{\frac{1}{2}\int_{0}^{T} \frac{1}{\s^2(Y_{\eta(t)\land\cdot}^i, \mu_{\eta(t)}^M) }\left(b(\t, Y_{\eta(t)\land\cdot}^i, \mu_{\eta(t)}^M) - b(\t_0, Y_{\eta(t)\land\cdot}^i, \mu_{\eta(t)}^M) \right)^2\dif t}\right\}.
\de
Thus, the maximum likelihood estimator of the parameter $\t$ is given by
\be
\t_T^M := \arg\max_{\t \in \Theta}{L^M_{T}(\t)}.\label{52} 
\ee

\section{An example}\label{exam}

In the section, we present an example to explain our results.

Consider the following MVSDE on $\mR$:
\be
\dif X_t = (\t X_t + \b \mE[X_t]) \dif t + \s\dif W_t, \quad X_0 = x_0\in\mR,
\label{61} 
\ee
where $\t \in \Theta$ is a unknown parameter and $\b, \s$ are nonzero constants.
Using the numerical simulation method in Section \ref{numsim}, we have the following numerical equation for Eq.(\ref{61})
\be\left\{\begin{array}{l}
Y_{0}^{i}=x_{0},\\
{Y_t}^{i}=Y_{t_k}^{i}+\left(\t Y_{t_{k}}^{i}+\b \frac{1}{N}\sum_{j=1}^{N} Y_{t_{k}}^{j}\right)(t-t_k) +\s(W^i_t-W_{t_k}^i) , \quad t\in[t_k, t_{k+1}].
\end{array}
\right.
\label{numsimeq}
\ee
In terms of the number $N$ of particles and the step size $M$, we draw Figure 1 and Figure 2. That is, we take $N=160, M=16$ in Figure 1, and $N=2560, M=256$ in Figure 2. Comparing Figure 1 with Figure 2, one can find that the numerical solution has higher frequency and smaller amplitude when the number of particles and the step size are larger.

\begin{figure}[ht]
\centering
\includegraphics[scale=0.4]{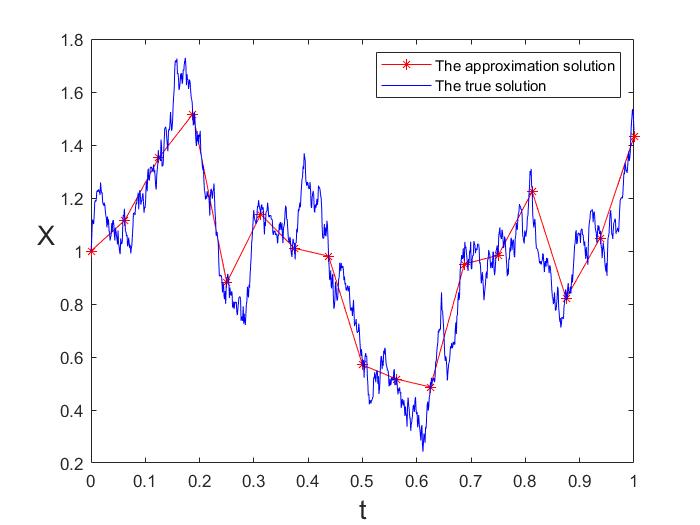}
\caption{Comparison of approximate solution and true solution, taking $\b=\sigma=1, N=160, M=16$.}
\label{t16}
\end{figure}

\begin{figure}[ht]
\centering
\includegraphics[scale=0.4]{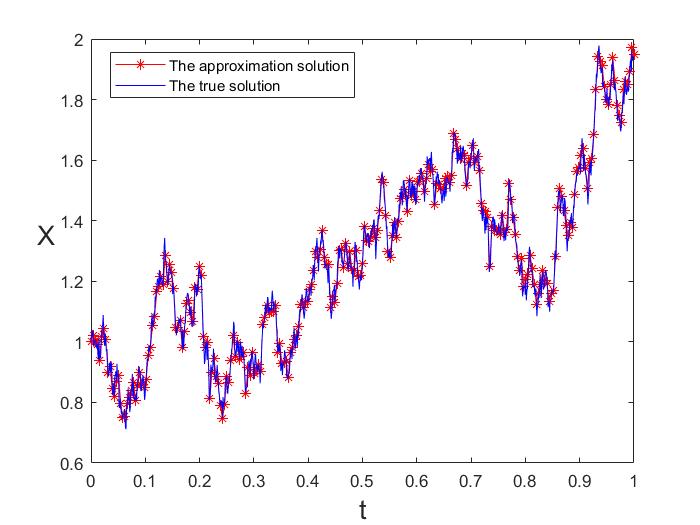}
\caption{Comparison of approximate solution and true solution, taking $\b=\sigma=1, N=2560, M=256$.}
\label{t256}
\end{figure}

According to (\ref{errest}) in Section \ref{numsim}, we calculate the errors between the solutions of Eq.(\ref{numsimeq}) and the solution of Eq.(\ref{61})  and list them in Table 1. From Table 1, one can find that the error decreases when the number of particles and the step size increase. 
\begin{table}[ht]
    \centering
    \begin{center} 
    \fontsize{12}{12}\selectfont    
    \caption{The errors between the numerical solution and the original solution when N, M take different values.}
    \setlength{\tabcolsep}{4mm}
     \begin{tabular}{ cccccc }
   \toprule
\diagbox [width=1.8cm, trim=l] {$N$}{$M$} & 16 & 32 & 64 & 128 & 256 \\
   \midrule  
		 160&0.0753     &0.0389    &0.0182   &0.0093    &0.0040\\
    		 320&0.0686    &0.0354    &0.0176    &0.0086    &0.0048\\
   		 640&0.0656    &0.0337    &0.0167    &0.0077    &0.0037\\
   		1280&0.0670    &0.0331     &0.0158    &0.0077    &0.0034\\
    		2560&0.0672    &0.0323    &0.0157    &0.0073    &0.0032\\
    \bottomrule
\end{tabular}
\label{error}
\end{center}
\end{table}

Finally, by the formula (\ref{52}) in Section \ref{numsim}, we get the maximum likelihood estimator $\t^M_T$ as follows:
\be
\t^M_T = \frac{\sum\limits_{k=0}^{M-1} Y_{t_{k}}^{i}(Y_{t_{k+1}}^{i}-Y_{t_k}^{i})-\sum\limits_{k=0}^{M-1} \b Y_{t_{k}}^{i}·\frac{1}{N}\sum\limits_{j=1}^{N} Y_{t_{k}}^{j}\frac{T}{M}}{\sum\limits_{k=0}^{M-1} (Y_{t_{k}}^{i})^2\frac{T}{M}}.\label{63} 
\ee
In terms of $T$, the values of $\t^M_T$ present in Table 2, which indicates that the value of $\t_T$ is closer to the true value $\t_0 = -0.5$ when the time $T$ becomes larger.
\begin{table}[ht]
    \centering
    \fontsize{12}{12}\selectfont    
    \caption{The maximum likelihood estimator $\t^M_T$ with $\b=\s=1, N=2560, M=256$.}
    \setlength{\tabcolsep}{4mm}
    \begin{tabular}{ cccccc }
    \toprule
   		$T$ & 1 & 2 & 5 & 8 & 10 \\
   \midrule  
	    $\t^M_T$ &-1.0510  &-0.7420   &-0.5107   &-0.5009   &-0.4999\\
    \bottomrule
\end{tabular}
\label{esti}
\end{table}

\bigskip

\textbf{Acknowledgements:}

The authors are very grateful to Professor Xicheng Zhang for valuable discussions. The second author also thanks Professor Renming Song for providing her an excellent environment to work in the University of Illinois at Urbana-Champaign.


\begin{thebibliography}{999}

\bibitem{ay} Y. Ait-Sahalia: Nonparametric pricing of interest rate derivative securities. {\it Econometrica}, 64(1996)527-560.

\bibitem{bjpn} J. P. N. Bishwal: Large deviations inequalities for the maximum likelihood estimator and the Bayes estimators in nonlinear stochastic differential equations. {\it Statistics} \& {\it Probability Letters}, 43(1999)207-215.

\bibitem{boss} M. Bossy, D. Talay: A stochastic particle method for the McKean-Vlasov and the Burgers equation. {\it Mathematics of Computation}, 66(1997)157-192.

\bibitem{cf} R. Carmona, F. Delarue: {\it Probabilistic Theory of Mean Field Games with Applications I}. Springer, 2018.

\bibitem{d} D. Dawson, J. Vaillancourt: Stochastic McKean-Vlasov equations. {\it Nonlinear Differential Equations and Applications NoDEA}, 2(1995)199-229.

\bibitem{dq1} X. Ding, H. Qiao: Euler-Maruyama approximations for stochastic McKean-Vlasov equations with non-Lipschitz coefficients, https://arxiv.org/abs/1903.11754.

\bibitem{dq2} X. Ding, H. Qiao: Stability for stochastic McKean-Vlasov equations with non-Lipschitz coefficients，https://arxiv.org/abs/1905.07883.

\bibitem{lipster} R. S. Liptser, A. N. Shiryayev: {\it Statistics of Random Processes: I General Theory}. 2nd Edition, Springer, 2001.

\bibitem{mckean} H. P. McKean: A Class of Markov Processes Associated with Nonlinear Parabolic Equations. {\it Proceedings of
the National Academy of Sciences}, 56(1966)1907-1911.

\bibitem{q} H. J. Qiao: A Nonlinear Stochastic Evolution Equation in Hilbert space. {\it Acta Mathematica Scientia}, 29A(2009)383-391.

\bibitem{r} P. Ren, J. Wu: Least squares estimator for path-dependent McKean-Vlasov SDEs via discrete-time observations. {\it Acta Mathematica Scientia}, 39B(2019)691-716.

\bibitem{s} B. M. Bibby, M. S$\phi$rensen: Martingale estimation functions for discretely observed diffusion processes. {\it Bernoulli}, 1(1995)17-39.

\bibitem{w} J. Wen, X. Wang, S. Mao and X. Xiao: Maximum likelihood estimation of McKean-Vlasov stochastic differential equation and its application. {\it Applied Mathematics and Computation}, 274(2016)237-246.

\bibitem{y} N. Yoshida: Estimation for diffusion processes from discrete observation. {\it Journal of Multivariate Analysis}, 41(1992)220-242.

\bibitem{zx} X. Zhang: Euler-Maruyama approximations for SDEs with non-Lipschitz coefficients and applications, {\it J. Math. Anal. Appl.,} 2006, 316: 447-458.
\end{thebibliography}
\end{document}